\documentclass[11pt]{article}

\usepackage{latexsym,amsmath,amsfonts,graphicx,graphics,epsfig}
\usepackage{theorem}
\newcommand{\nref}[1]{(\ref{#1})}
\def\Cal{\cal}

\def\text#1{\hbox{#1}}

\def\endproof{\mbox{\ $\Box$}}

\newcommand{\R}{\mathbb R}

\def\1{\mbox{1\hspace{-.20em}I}}

\newtheorem{theorem}{Theorem}[section]
\newtheorem{proposition}{Proposition}[section]
\newtheorem{lemma}{Lemma}[section]

\newtheorem{remark}{Remark}[section]

\newcommand{\CC}{{\Cal{C}}}

\newcommand{\CS}{{\Cal{S}}}

\def\d_1{\gamma_1}
\def\l{\left}
\def\r{\right}

\voffset=-20mm \numberwithin{equation}{section}

\begin{document}
\title{Sharp Variable Selection of a Sparse Submatrix in a High-Dimensional Noisy Matrix}

\author{Butucea, C. \thanks{Universit\'e Paris-Est, LAMA (UMR 8050), UPEMLV, UPEC, CNRS, F-77454, Marne-la-Vall\'ee, France  and CREST, Timbre J340 3,
av.\ Pierre Larousse, 92240 Malakoff Cedex, France}, \boxed{ \rm{Ingster, Yu. I.}}
 and Suslina, I. A.
\thanks{ St.Petersburg  National Research University of Information Technologies, Mechanics and Optics,
 49 Kronverkskiy pr., 197101, St.Petersburg, Russia} }


\maketitle

\begin{abstract}

We observe a $N\times M$ matrix of independent, identically distributed
Gaussian random variables which are centered except for elements of some submatrix of size $n\times m$ where the mean is larger than some $a>0$. The submatrix is sparse in the sense that $n/N$ and $m/M$ tend to 0, whereas $n,\, m, \, N$ and $M$ tend to infinity.

We consider the problem of selecting the random variables with significantly large mean values. We give sufficient conditions on $a$ as a function of $n,\, m,\,N$ and $M$ and construct a uniformly consistent procedure in order to do sharp variable selection.
We also prove the minimax lower bounds under necessary conditions which are complementary to the previous conditions. The critical values $a^*$ separating the necessary and sufficient conditions are sharp (we show exact constants).

We note a gap between the critical values $a^*$ for selection of variables and that of detecting that such a submatrix exists given by \cite{ButIn}. When $a^*$ is in this gap, consistent detection is possible but no consistent selector of the corresponding variables can be found.

\end{abstract}

\medskip
\noindent {\bf Keywords:} estimation, minimax testing, random matrices, selection of sparse signal, sharp selection bounds, variable selection.

\newpage
\section{Introduction}

We observe random variables that form an $N\times M$ matrix
$\mathbf{Y}=\{Y_{ij}\}_{i=1,...,N, j=1,...,M}$:
\begin{equation}\label{mod}
Y_{ij}=s_{ij}+  \xi_{ij},\quad i=1,\ldots,N,\quad j=1,\ldots,M,
\end{equation}
where  $\{\xi_{ij}\}$ are i.i.d. random variables and $s_{ij}\in\R$,
for all $i\in\{1,...,N\}$, $j\in \{1,...,M\}$. The error terms
$\xi_{ij}$ are assumed to be distributed as standard Gaussian $\mathcal{N}(0,1)$ random variables.

Let us denote by
\begin{equation}\label{subset}
\mathcal{C}_{nm} = \left\{C=A\times B \subset \{1,\ldots,N\} \times \{1,\ldots,M\},\, Card(A)=n,\, Card(B)=m  \right\},
\end{equation}
the collection of subsets of $n$ rows and $m$ columns of a matrix of size
$N \times M$.

We assume that our data have mean $s_{ij} = 0$ except for elements in a submatrix of size $n \times m$, indexed by a set $C_0$ in $ \mathcal{C}_{nm}$, where $s_{ij} \geq a$, for some $a>0$.


Our model means that, for some
$a>0$ which may depend on $n,\, m, \, N$ and $M$,
\begin{equation}\label{bin}
 \exists\ C_0\in{\cal{C}}_{nm}\ \text{such that}\ \ s_{ij}=0,\
\text{if}\ (i,j)\notin C_0,\ \text{and}\ s_{ij} \geq a,\ \text{if}\
(i,j)\in C_0.
\end{equation}

Let ${{\cal S}}_{nm,a}$ be the collection of all matrices $S=S_C$, $C \in \mathcal{C}_{nm}$ that satisfy \nref{bin}. Our model implies also that there exists some $C_0$ in $\mathcal{C}_{nm}$ such that $S=S_{C_0}$ belongs to $\mathcal{S}_{nm,a}$.

We discuss here only significantly positive means of our random variables. The problem of selecting the variables with significantly negative means can be treated in the same way, by replacing variables $Y_{ij}$ with $ - Y_{ij}$.

\bigskip

Denote by $P_C$ the
probability measure that corresponds to observations \nref{mod} with
matrix $S=S_C = \{s_{ij}\}_{i=1,...,N,\, j=1,...,M}$, $s_{ij}=0$ if $(i,j) \not \in C$, $s_{ij} \geq a >0$ if $(i,j) \in C$. We also denote $P_0 = P_{C_0}$ and $E_0$ the expected value with respect
to the measure $P_0$.

Our goal is to propose a consistent estimator of $C_0$, that is to select the variables in the large matrix of size $N \times M$ where the mean values are significantly positive. Our approach is to find the boundary values of $a>0$, as function of $n, \, m, \, N$ and $M$, where consistent selection is possible and separate them from the cases where consistent selection is not possible anymore.

We are interested here in sparse matrices, i.e. the case when $n$ is
much smaller than $N$ and $m$ is much smaller than $M$.

\bigskip

Large data sets of random variables appear nowadays in many applied fields such as signal processing, biology and, in particular, genomics, finance etc. In genomic studies of cancer we may require to detect sample-variable associations see~\cite{SWPN}. Our problem further adresses the question: if such an association is detected can we estimate the sample components and the particular variables involved in this association?

We may also view our problem as a matrix-mixture model, where each observation $Y_{ij}$ has distribution
$$
Y_{ij} \sim ( 1-p)\cdot \mathcal{N}(0,1) + p\cdot \mathcal{N}(s_{ij},1),
$$
with $p=p_{n,N,m,M} \in (0,1) $ the mixture probability (small) and $s_{ij} \geq a$ for $(i,j) \in C_0$.
Such models appear, for example, in multiple testing setup where $Y_{ij}$ are test statistics, which are i.i.d. under the null hypothesis and they have a Gaussian distribution. Benjamini and Hochberg~\cite{BH} proposed to study the false discovery rate and many models have been proposed since for estimating $p$ and the mixture density of the observations in the non Gaussian case.
In our approach the multiple tests are indexed by $(i,j) \in \{1,...,N\}\times\{1,...,M\}$ such that the mixture occurs with a submatrix structure.
We address here the question of selecting the multiple tests which are significant (have rejected the null) in a matrix setup, and, as a particular case, in a vector setup as well. This problem is also known as classification and
it was known that in some cases classification is not possible even though detection is possible, see \cite{CJL}. Our result provides new rates for the matrix case and sharp constants for the vector case.

\bigskip

Sparsity assumptions were introduced for vectors. There is a huge amount of literature for variable selection in (sparse or not sparse) linear and nonparametric regression, gaussian white noise and density models. Estimation of the sparse vector as well
as hypothesis testing for vectors were thoroughly studied under various sparsity assumptions as well. See for example Bickel, Ritov and Tsybakov~\cite{BRT}
and references therein, for estimation issues, and Donoho and Jin~\cite{DJ04}, Ingster~\cite{I97} and Ingster and Suslina~\cite{IS02b}, for testing.

In the context of matrices, different sparsity assumptions can be
imagined. For example, matrix completion for low rank matrices with
the nuclear norm penalization has been studied by Koltchinskii,
Lounici and Tsybakov~\cite{KLT}.

The detection problem was considered in this setting by
Butucea and Ingster ~\cite{ButIn}. A more general setup, where each observation is replaced by a smooth signal was considered by Butucea and Gayraud~\cite{BG}. We can apply our results to their setup in order to select the signals with significant energy (norm larger than $a$).

\bigskip

We study here the variable selection problem in a matrix from a minimax point of view. A selector is
any measurable function of the observations, $\hat C=
\hat C(\{Y_{ij}\})$ taking values in ${\cal{C}}_{nm}$. For such a
selector $\hat C=\Hat C(Y),\ Y=\{Y_{ij}\}$ we denote the maximal risk by
$$
R_{nm,a}(\hat C)=\sup_{{S_{C_0}\in\,\CS_{nm,\,a}}}P_{C_0}(\hat C(Y)\ne C_0).
$$

We define the minimax risk  as
$$
R_{nm,a}=\inf_{\hat C}R_{nm,a}(\hat C).
$$

From now on, we assume in the asymptotics that $N\to\infty,\
M\to\infty$ and $n=n_{NM}\to\infty,\ n\ll N, \ m=m_{NM}\to\infty,\
m\ll M $. Other assumptions will be given later.

We say that a selector is consistent in the minimax sense, if $
R_{nm,a}(\hat C) \to 0 $.
\bigskip

We suppose that $a>0$ is unknown. The aim of this paper is to give
asymptotically sharp boundaries for minimax selecting risk. It means
that, first, we are interested in the conditions on $a=a_{NM}$ which
guarantee the possibility of selection i.e., the fact that $R_{nm,a}
\to 0 $. We construct the selecting procedure
\begin{equation}\label{Y}
\hat C^{\star}(Y)=\arg\max_{C\in\,{\cal{C}}_{nm}}\sum_{(i,j)\in C}Y_{ij}
\end{equation}
We investigate the upper bounds of the minimax selection risk of
this procedure. Second, we describe conditions on $a$ for which we
have the impossibility of selection, i.e., the lower bounds
$R_{nm,a}\to 1 $. These results are called the lower bounds. The two
sets of condition are partially complementary in a sense that
violation of the upper bound conditions
 imply either impossibility of selection or indistinguishability
(see ~\cite{ButIn}).
\begin{remark}
Note that $P_{C_0}(\hat C^\star(Y)\ne C_0)$  does not depend on
$C_0=C_0(N,M,n,m,a)$. Therefore, for any $C_0$ we have
$$
R_{nm,a}(\hat C^\star) = \max_{{S_{C_0}\in\,\CS_{nm,\,a}}}
P_{C_0}(\hat C^\star(Y) \neq C_{0})=P_{C_0}(\hat C^\star(Y) \neq
C_{0}).
$$
\end{remark}

\bigskip

The problem of choosing a submatrix in a Gaussian random matrix has
been previously studied by Sun and Nobel~\cite{SN}. They are interested in the largest squared submatrix in $Y$ under the null hypothesis such that its average
 is larger than some fixed threshold.
The algorithm of choosing such submatrices was previously introduced
in Shabalin {\it et al.}~\cite{SWPN}.

\bigskip

The plan of the paper is as follows.
In Section~\ref{mres} we state the main results of this paper: the upper bounds for the selection procedure $\hat C^\star$ under conditions on $a$,
 as well as inconsistency property of this procedure under complementary conditions on $a$, and, finally, lower bounds for variable selection.
  We compare these results with the results for detection in~\cite{ButIn}. We give results for the vector case ($m=M=1$) which are new as far as the asymptotic
   constant is concerned.
 In Section~\ref{upb} we prove the upper bounds for the selection of variables, that is a bound from above on $a$, in which
  $R_{nm,a}(\hat C^\star) = \sup_{S_{C_0}} P_{C_0}(\hat C^\star \neq C_{0})\to 0$.
 In Section~\ref{lowb1} we prove lower bounds for variable selection, that is,
 a bound on the parameter $a$ from below which imply that the minimax estimation risk $R_{nm,a}$ tends to 1.
Two techniques provide the sharp lower bounds. One method is classical for nonparametric estimation, while the other makes a generalization of a well-known result
to testing $L\geq 2$ hypotheses: the minimax risk is larger than the risk of the maximum likelihood estimator.

Future extensions of this problem include several open problems. For example, consider two-sided variable selection,
 i.e. finding $C_0$ where the mean $|s_{ij}|\geq a$, for $(i,j) \in C_0$. Another possibility is to consider non Gaussian observations,
  but having distribution in the exponential family. As mentioned, we may replace each observation with a smooth signal and detect the active
   components (signals with significant total energy) in the matrix.

\section{Main Results}\label{mres}
Let
\begin{equation}\label{lim}
N\to\infty,\, n\to\infty,\, p=n/N\to 0;\quad M\to\infty,\, m\to\infty,\,
q=m/M\to 0.
\end{equation}


We suppose that $a>0$ is unknown. The aim of this paper is to give
asymptotically sharp boundaries for variable selection in a sparse high-dimensional matrix.
Our approach is to give, on the one hand, sufficient asymptotic conditions on $a$ such that the probability
of wrongly selecting the variables in $C_0$ tends to 0 and, on the other hand, conditions
under which no consistent selection is possible.

First, we are interested in the conditions on $a=a_{nmNM}$ which
guarantee consistent variable selection, i.e., the fact that
we construct the selector $\hat C^\star$ in \nref{Y} and prove that
$R_{nm,a}(\hat C^\star) \to 0$.
The selector $\hat C^\star$ is scanning the large $N\times M$ matrix and maximizes the sum
of the inputs over all $n \times m$ submatrices.

The key quantities appearing in next theorems are
\begin{equation}\label{abc}
\begin{array}{ll}
B = B_{n,m,N,M} = \min\{A_{1}, \, A_2,\, A\}, &\mbox{ where }
A = \frac{a \sqrt{nm}}{ \sqrt{2(n \log(p^{-1})+m \log(q^{-1} ) )} },\\
& \\
A_1 = \frac{a \sqrt{m}}{ \sqrt{2} (\sqrt{\log(n)}+\sqrt{\log(N-n)}) },
& A_2 = \frac{a \sqrt{n}}{ \sqrt{2} (\sqrt{\log(m)}+\sqrt{\log(M-m)})}.
\end{array}
\end{equation}
Let us consider the particular case where the matrix and the submatrix are squared ($N=M$ and $n=m$) and, moreover, such that
$$
\frac{\log(n)}{\log(N)} = \frac{\log(m)}{\log(M)} \to 0.
$$
Then, $\log(n(N-n))\sim \log(N/n)$ and $\log(m(M-m)) \sim\log(M/m)$
which imply that $A_1=A_2 \geq A$ and, therefore, $B=A$. We need
terms $A_1=A_2$ in order to consider cases where $\liminf
\log(n)/\log(N) $ and $\liminf \log(m)/\log(M) $ are large enough and close to 1.

Another particular example is
$n\sim N^P$ or $m\sim M^Q$, for $P,\,Q \in (0,1)$ that we discuss in more details later on.

For this reason, we distinguish the case of severe sparsity when $B=A$, from the case of moderate sparsity when $B=A_1$ or $B=A_2$.

\medskip

The following Theorem gives sufficient conditions for the
 boundary $a=a_{n,m,N,M}$ such that selection is consistent uniformly over tha class $\mathcal{S}_{nm,a}$. The
selector which attains these bounds is defined by \nref{Y}.

\begin{theorem}\label{ub} {\bf Upper bounds. }
 Assume \nref{lim} and assume $B = B_{n,m,N,M}$ defined by \nref{abc} is such that
 \begin{equation}\label{B}
 \liminf B_{n,m,N,M} > 1,
 \end{equation}
 then the selector $\hat C^\star$ given by \nref{Y} is
consistent, that is
$$
R_{nm,a}(\hat C^\star)=P_{C_0}(\hat C^\star \ne C_0)\rightarrow 0.
$$
\end{theorem}
{\bf Proof} is given in Section \ref{upb}.

Condition \nref{B} is equivalent to saying that
$$
\liminf A>1 \mbox{ and } \liminf A_1 >1  \mbox{ and } \liminf A_2 >1.
$$
The following proposition says that $\liminf A_1 >1$ and $\liminf A_2 >1$ are necessary conditions for the consistency (in the minimax sense) of
 the selector $\hat C^\star$ of $C_0$.
\begin{proposition} \label{notcons}Assume \nref{lim} and let the selector $\hat C^\star$ be the selector given by \nref{Y}.
If
$$
\limsup A_1 < 1
\text{ or }
\limsup A_2 < 1
$$
then, for any $C_0$ such that $S_{C_0} \in \mathcal{S}_{nm,a}$,
$$
P_{C_0}(\hat C^\star \not = C_0) \to 1.
$$
\end{proposition}
{\bf Proof} is given in Section~\ref{sec:notcons}.

In the following theorem we give a sufficient condition on $a$ under which consistent selection of $C_0$ is impossible uniformly over
 the set $\mathcal{S}_{nm,a}$. These are the minimax lower bounds for variable selection.

\begin{theorem} \label{estimlb}
Assume \nref{lim}. If, moreover, $B = B_{n,m,N,M}$ defined by \nref{abc} is such that
\begin{equation} \label{Bsup}
\limsup B_{n,m,N,M} < 1,
\end{equation}
then there is no consistent selection of $C_0$ uniformly over $\mathcal{S}_{nm,a}$, that is
$$
\inf_{\hat C} \sup_{{S_{{C_0}}\in\,\CS_{nm,\,a}}}P_{{C_0}}(\hat C(Y)\ne {C_0}) \to 1,
$$
asymptotically, where the infimum is taken over all measurable functions $\hat C = \hat C(Y)$.
\end{theorem}
{\bf Proof} of this theorem is given in Section~\ref{sec:estimlb} and \ref{sec:notcons}.

Theorems \ref{ub} and \ref{estimlb} imply that the critical value for $a$ is
\begin{eqnarray}
a^* &\sim &\max\left\{\frac{\sqrt{2\log(n)}+\sqrt{2 \log(N-n)}}{\sqrt{m}},
\frac{\sqrt{2\log(m) }+\sqrt{2\log(M-m)}}{\sqrt{n}},\right. \nonumber\\
&& \left.
\frac{\sqrt{2(n\log(N/n) + m \log(M/m))}}{\sqrt{nm}} \right\}.\label{astar}
\end{eqnarray}
By critical we mean in the sense that, for $a$ such that $\lim\inf a/a^\star >1$, there is an estimator which is uniformly consistent, while, for $a$
 such that $\lim\sup a/a^\star <1$, no uniformly consistent estimator exists.

If we consider the particular case where $n = N^{P}$ and $m  = M^{Q}$ grow polynomially, for some fixed $P, \, Q$ in $(0,1)$, the critical value becomes
\begin{eqnarray*}
(a^*)^2 &\sim &\max\left\{ \frac{2(1+\sqrt{P})^2\log(N)}{m}, \frac{2(1+\sqrt{Q})^2\log(M)}{n} , \right.\\
&& \left. \frac{2(1-P)\log(N)}m + \frac{2(1-Q) \log(M)}n
\right\} .
\end{eqnarray*}
 If, moreover, $n=m$ and $N=M$, we get $(a^*)^2 \sim \max\{ {2} (1+\sqrt{P})^2, 4(1-P)\} {\log(N)/n}$. So,
 the amount of sparsity depends on whether $P$ is larger or smaller than 1/9. In this particular example, we have moderate sparsity, $B=A_1=A_2 \leq A$,
 as soon as $P \geq 1/9$.

\subsection{Variable selection vs. detection}

Let us compare the result in Theorem~\ref{ub} and
Theorem~\ref{estimlb} with the upper bounds and the lower bounds for
detection of a set $C_0$ where our observations have significant
means, i.e. above threshold $a$. The testing problem for our model
can be stated as
$$
H_0: s_{ij}=0 \mbox{ for all } (i,j)
$$
and we call $P_0$ the likelihood in this case, against the
alternative
$$
H_1: \mbox{ there exists } \, C_0 \in \CC_{nm} \mbox{ such that }
S=S_{C_0} \in \CS_{nm,a}.
$$

Recall the following theorems.

\begin{theorem}{\bf Upper bounds for detection}\label{teorupb}, see \cite{ButIn}.
  Assume \nref{lim} and let $a$ be such that at least one of
 the following conditions hold
\begin{equation*}\label{cond1}
a^2nmpq =\frac{(anm)^2}{NM}\to\infty
\quad \mbox{ or }\quad
\liminf A > 1.
\end{equation*}
Then distinguishability is possible, i.e.
$$
\inf_{\psi(Y)} \left( P_0( \psi(Y) = 1)+ \sup_{S_{C_0} \in
\CS_{nm,a}} P_{C_0}(\psi(Y) = 0)\right) \to 0,
$$
where the infimum is taken over all measurable functions $\psi$
taking values in $\{0,1\}$.
\end{theorem}

%
%
%

It was also shown in \cite{ButIn}, that the asymptotically optimal
test procedure $\psi^*$ combines the scan statistic based on our
$\hat C^\star$ with a linear statistic which sums all observations
$Y=\{Y_{ij}\}_{i,j}$.
The test procedure $\psi^*$ rejects the null hypothesis as soon as
either the linear or the scan test rejects.

\bigskip

\begin{theorem}\label{lb1} {\bf Lower bounds for detection}, see \cite{ButIn}.
 Assume \nref{lim} and
\begin{equation}\label{cond3a}
{n}\log(p^{-1})\asymp m\log(q^{-1}), \quad
\frac{\log\log(p^{-1})}{\log(q^{-1})}\to 0,\quad
\frac{\log\log(q^{-1})}{\log(p^{-1})}\to 0.
\end{equation}
Moreover, assume that
\begin{equation*}\label{linsup}
a^2 nmpq = \frac{(anm)^2}{NM} \to 0
\quad \mbox{ and } \quad
\limsup A <1.
\end{equation*}
Then, consistent detection is impossible, that is
$$
\inf_{\psi(Y)} \left( P_0( \psi(Y) = 1)+ \sup_{S_{C_0} \in
\CS_{nm,a}} P_{C_0}(\psi(Y) = 0)\right) \to 1,
$$
where the infimum is taken over all measurable functions $\psi$
taking values in $\{0,1\}$.
\end{theorem}

We deduce that there is a gap between least conditions for testing
that $C_0$ exists and selection of the actual variables $(i,j) \in
C_0$ (estimation of $C_0$). In Table~\ref{table} we summarize possible cases
were consistent selection and/or consistent testing is possible or not.
We can prove that, if
\begin{equation*}
 \limsup A <1, \quad \lim\inf A_1 >1 \quad  \mbox{ and } \quad \lim\inf A_2 >1
\end{equation*}
then $a^2 nmpq \to 0$, hence Theorem~\ref{lb1}.
We used this in the conditions of the second case where neither consistent selection,
nor testing is possible.

\begin{table}[hptb!]\label{table}
\begin{center}
\begin{tabular}{c|c|c}
Selection $\backslash$ Test& Yes & No\\
\hline
&&\\
Yes & $\liminf B >1$ & -\\
&&\\
\hline
No &
\begin{tabular}{l}
1) $\, \limsup B <1$ \\
and $ a^2nmpq \to \infty$\\
\\
2) $\, \liminf A >1$ and \\
  ($\limsup A_1 <1$ or $\limsup A_2 <1$)
\end{tabular}
&
\begin{tabular}{l}
Under \nref{cond3a} for the test: \\
1) $\,\limsup A <1 $ \\
and $ a^2nmpq \to 0 $\\
\\
2) $\, \limsup A <1$ and\\
$\, \liminf A_1 >1$ and \\
$\, \liminf A_2 >1$
\end{tabular}
\end{tabular}
\end {center}
\caption{Conditions for variable selection and/or testing}
\end{table}

%
%
%
%
%
%
%

%
%
%
%
%

\medskip

Let us consider the following example: $N=n^2,\ M=\log(n), \, \
m=\log\log(n) $ (and, for instance, $a^2=\log(n)/\log\log(n)$). For
all $a$ such that $a^2\gg\log(n)/(\log\log(n))^2$ as $n\to \infty$,
we have $a^2nmpq=a^2(\log\log(n))^2/\log(n)\to\infty$.
 Therefore, on the one hand,
distinguishability holds, see Theorem \nref{teorupb}, i.e. we can
construct a particular test procedure $\psi^\star$ such that
$$
P_0( \psi^\star(Y) = 1)+ \sup_{S_{C_0} \in
\CS_{nm,a}} P_{C_0}(\psi^\star(Y) = 0) \to 0.
$$
On the other hand,
$$
\frac{a^2m}{2(\sqrt{\log(n)}+\sqrt{\log(N-n)})^2}=\frac{a^2\log\log(n)}{(2+\sqrt{2})^2\log(n)}(1+o(1))<1,
$$
for all $a$ such that $a^2<(1-\delta)(2+\sqrt{2})^2
\log(n)/\log\log(n),\ \delta>0$. By Theorem 2.2, no consistent
selection is possible in this case.

\bigskip

\subsection{Vector case}

Previous results can also be proven for the vector case, that is for the gaussian independent, observations
$$
X_i = s_i +\xi_i,\quad i=1,...,N,
$$
where $s_i \geq a$ for all $i$ in a set $A_0$ of $n$ elements and $s_i=0$ otherwise. We suppose $n,\, N \to \infty$ such that $n/N \to 0$.
Similarly, we can show the following result.
\begin{theorem} {\bf Upper bounds} In the previous model, if
$$
\liminf \frac{a}{\sqrt{2 \log (N)}+\sqrt{2\log(n)}}>1,
$$
then the estimator $\hat A^\star = \arg\max_A \sum_{i\in A} X_i$ is such that
$$
\sup_{A_0} P_{A_0}(\hat A^\star \ne A_0) \to 0.
$$
{\bf Lower bounds} If
$$
\limsup \frac{a}{\sqrt{2 \log (N)}+\sqrt{2\log(n)}} < 1,
$$
then
$$
\inf_{\hat A}\sup_{A_0} P_{A_0}(\hat A \ne A_0) \to 1.
$$
\end{theorem}
The critical value is $a^\star=\sqrt{2 \log N}+\sqrt{2\log(n)}$. It
is equivalent to $\sqrt{2 \log N}$ if $\log(n)/\log(N)\to 0$ and
$a^\star = \sqrt{2}(1+\sqrt{1-\beta})\sqrt{ \log N}$ if
$N=n^{\beta}$ for some $\beta\in (0,1)$. This result
follows from \cite{ingst} (see
Section 3.1, Remark 2 and references therein).

Note that in the vector case, variable selection was mostly studied for the regression model with deterministic design,
see e.g. \cite{ACL}, \cite{Verzelen} and references therein.

Our results  are sharp as they give also the asymptotic constant.

Let us stress the fact that the particular case we study here is fundamentally different from the matrix setup.
 Indeed, an additional regime is observed according to the sparsity structure of the submatrix (severe or moderate) and it cannot be obtained from previous results for vectors by, say, vectorizing the matrix.

\section{Upper bounds}\label{upb}

{\bf Proof of Theorem \ref{ub}} Note that
$$
P_{C_0}(\hat C^\star \ne C_0)
=P_{C_0}(\max_{C\in\CC_{nm}}\sum_{(i,j) \in C} Y_{ij} - \sum_{(i,j) \in C_0} Y_{ij}>0).
$$
We shall split the sets $C$ according to the size of their common elements with the true underlying $C_0$. Let $C=A\times B$ and $C_0 = A_0 \times B_0$ and let $k$ be the number of elements in $A \cap A_0$ and $l$ the number of elements in $ B\cap B_0$. Then, if we denote by $\CC_{nm,kl}$ the collection of such matrices $C$:
\begin{eqnarray}
P_{C_0}(\hat C^\star \ne C_0)
& = & P_{C_0} \left( \max_{k=0,...,n} \max_{l=0,...,m} \max_{C \in\, \CC_{nm,kl}}  \sum_{(i,j) \in C} Y_{ij} - \sum_{(i,j) \in C_0} Y_{ij}>0\right) \nonumber \\
& \leq &   P_{C_0} \left(\max_{k=0,...,n} \max_{l=0,...,m}
\max_{C \in\, \CC_{nm,kl}}(\sum_{C\setminus  C_0} \xi_{ij} - \sum_{C_0 \setminus C} \xi_{ij} - a(nm-kl)) > 0 \right).\nonumber
\end{eqnarray}
From now, we fix $0<\delta <1$ and separate two cases: when $kl <(1-\delta )nm$ and when $kl \geq (1-\delta) nm$. As $\delta$ will be chosen small, it means that we treat differently the cases where the matrix $C$ overlaps $C_0$ but weakly (or not at all) and where the matrices overlap almost entirely.
We write and deal successively with each term in
\begin{eqnarray}
\label{weakdoublesum}
&&P_{C_0}(\hat C^\star \ne C_0) \nonumber \\
&\leq & P_{C_0} \left(\max_{k,l} \max_{kl<(1-\delta)nm}
\max_{C \in\, \CC_{nm,kl}}(\sum_{C\setminus  C_0} \xi_{ij} - \sum_{C_0 \setminus C} \xi_{ij} - a(nm-kl)) > 0 \right)\\
&&+
P_{C_0} \left(\max_{k,l} \max_{kl \geq (1-\delta)nm}
\max_{C \in\, \CC_{nm,kl}}(\sum_{C\setminus  C_0} \xi_{ij} - \sum_{C_0 \setminus C} \xi_{ij} - a(nm-kl)) >  0 \right).\label{strongdoublesum}
\end{eqnarray}

\subsection{Weak intersection}

Let us fix $k$ and $l$ such that $kl <(1-\delta) nm$ for some $0<\delta <1$.
Equivalently, we have $nm-kl > \delta nm$.
In this case, we shall bound the probability in \nref{weakdoublesum} as follows
\begin{eqnarray*}
&& P_{C_0} \left(\max_{k,\,l} \max_{kl<(1-\delta)nm}
\max_{C \in\, \CC_{nm,kl}}(\sum_{C\setminus  C_0} \xi_{ij} - \sum_{C_0 \setminus C} \xi_{ij} - a(nm-kl))>0 \right)\\
&\leq & \sum_{k=0}^n \sum_{l=0}^m I_{kl<(1-\delta)nm}P_{C_0} \left(
\max_{C \in\, \CC_{nm,kl}}\sum_{C\setminus  C_0} \xi_{ij}
+\max_{C \in\, \CC_{nm,kl}} \sum_{C\cap  C_0} \xi_{ij}
-\sum_{C_0} \xi_{ij} \geq a(nm-kl)
\right)\\
&\leq & \sum_{k=0}^n \sum_{l=0}^m I_{kl<(1-\delta)nm}(T_{1,kl}+T_{2,kl}+T_{3,kl}) ,
\end{eqnarray*}
where we denote by $I_{kl<(1-\delta)nm}$ the indicator function of the set where $kl<(1-\delta)nm$ and by
\begin{eqnarray*}
T_{1,kl} & = & P_{C_0} \left(
\max_{C \in\, \CC_{nm,kl}}\sum_{C\setminus  C_0} \xi_{ij} >(1-\delta_1)a(nm-kl)
\right)\\
T_{2,kl}&=& P_{C_0} \left(
\max_{C \in\, \CC_{nm,kl}}\sum_{C\cap C_0} \xi_{ij} >\frac{\delta_1}2 a(nm-kl)
\right)\\
T_{3,kl}&=& P_{C_0}\left(-\sum_{C_0} \xi_{ij} > \frac{\delta_1}2 a(nm-kl) \right),
\end{eqnarray*}
for some $0<\delta_1<1$.

Before continuing the proof, recall that, if $n,\, N$ tend to infinity, such that $n/N \to 0$, we have
\begin{eqnarray*}
\log(C_{N-n}^{n-k})&\sim &(n-k)\log \left(\frac{N-n}{n-k}\right) + (N-2n+k)\log \left(\frac{N-n}{N-2n+k} \right)\\
&\sim & (n-k)\log \left( \frac{N-n}{n-k} \right)
\end{eqnarray*}
and
$$
\log(C_n^k) \leq \min \left\{(n-k ) \log\left( \frac{ne}{n-k} \right) ,\, k \log \left( \frac{ne}k \right) \right\},
$$
for all $k=1,...,n-1$ and $\log C_n^n = 0$.

\bigskip

In order to give an upper bound for $T_{1,kl}$, we shall distinguish the case where $k< (1-\delta)n$ and $l=m$ (the case $k=n$ and $l < (1-\delta)m$ is treated similarly) from the case $kl < (1-\delta)nm$, $k<n$ and $l<m$.
On the one hand, if $k< (1-\delta)n$ and $l=m$, we write, for a generic standard gaussian random variable $Z$ (which might change later on):
\begin{eqnarray*}
T_{1,km} &\leq & P_{C_0} \left( \max_{A \in\,
\CC_{n,k}}\sum_{(A\setminus A_0) \times B_0} \xi_{ij}
>(1-\delta_1)a(n-k)m
\right)\\
&\leq & C_{N-n}^{n-k} P(Z >(1-\delta_1) a\sqrt{(n-k)m})\\
& \leq & \exp\left( -\frac {(1-\delta_1)^2}2 a^2 (n-k)m +
\log(C_{N-n}^{n-k}) \right),
\end{eqnarray*}
where we use repeatedly that $P(Z >u) \leq \exp(-u^2/2)$, for all $u
\geq 0$.

 Now,
\begin{eqnarray*}
\log(C_{N-n}^{n-k} ) &\leq & (n-k) \log \left(\frac{N-n}{n-k}
\right)(1+o(1)).
\end{eqnarray*}
Therefore,
\begin{eqnarray*}
T_{1,km} &\leq & \exp\left( -(n-k)\left(\frac {(1-\delta_1)^2}2 a^2
m  - \log\left(\frac{N-n}{n-k}\right)(1+o(1)) \right)\right).
\end{eqnarray*}

By assumption \nref{B} we can say that
\begin{eqnarray} \label{alpha}
&&\min\left\{\frac{a^2nm}{2(n\log(p^{-1})+m \log(q^{-1}))}, \right.
\nonumber \\
&&\left.
\frac{a^2m}{2(\sqrt{\log(N-n)}+\sqrt{\log(n)})^2},\frac{a^2n}{2
(\sqrt{\log(M-m)}+\sqrt{\log(m)})^2} \right\}\geq 1+\alpha,
\end{eqnarray}
for some fixed small $\alpha > 0$. Therefore, if $\delta_1>0$ is
small enough, we have some $\alpha_1>0$ such that
\begin{eqnarray} \label{alpha1}
\frac{(1-\delta_1)^2}2 a^2 m \geq (1+\alpha_1) ( \log(( N-n)n)
)>\log\left(\frac{N-n}{n-k}\right)(1+o(1))+\log(n),
\end{eqnarray}
asymptotically. Indeed, it is sufficient that $(1-\delta_1)^2
(1+\alpha) \geq 1+\alpha_1$.

We get
$$
T_{1,km} \le \exp(-(n-k)\log(n)) .
$$
We conclude that
$$
\sum_{k:(n-k)>\delta n} T_{1,km} \leq  n \max_{k:(n-k)>\delta
n}\{\exp(-(n-k)\log(n)\}<n^{-\delta n+1} =o(1).
$$

\bigskip

On the other hand, if $kl < (1-\delta)nm$, $k<n$ and $l<m$, note first that the maximum is taken over all $C$
 in $\CC_{nm,kl}$, but only the lines and columns outside $C_0$ actually play a role over the sum $\sum_{C\setminus C_0}\xi_{ij}$. There are
$C_{N-n}^{n-k} \cdot C_{M-m}^{m-l} \cdot C_n^k \cdot C_m^l$ different values of this sum. We write:
\begin{eqnarray}
T_{1,kl} & \leq & C_{N-n}^{n-k} \cdot C_{M-m}^{m-l} \cdot C_n^k \cdot C_m^l P\left(Z > (1-\delta_1) a\sqrt{nm-kl} \right) \nonumber\\
&\leq & C_{N-n}^{n-k} \cdot C_{M-m}^{m-l} \cdot C_n^k \cdot C_m^l\,
\exp\left( -\frac {(1-\delta_1)^2}2 a^2 (nm-kl) \right) \nonumber\\
&\leq & \exp\left(-\frac {(1-\delta_1)^2}2 a^2 (nm-kl)  + \log
(C_{N-n}^{n-k}C_{M-m}^{m-l} C_n^k  C_m^l) \right) \label{eqT1}.
\end{eqnarray}

As we have $n,\, m,\, N,\, M$ tend to infinity, then
\begin{eqnarray*}
\log(C_{N-n}^{n-k} \cdot C_{M-m}^{m-l}\cdot C_n^k \cdot C_m^l)
&\leq & \left((n-k) \log \frac{N-n}{n-k}+ (m-l) \log \frac{M-m}{m-l}
\right)(1+o(1))\\
&& + (n-k) \log \frac{ne}{n-k} + (m-l) \log \frac{me}{m-l}\\
&\leq & \left((n-k) \log \frac{N-n}{n}+ (m-l) \log \frac{M-m}{m}
\right)(1+o(1))\\
&& +\l((n-k) \log \frac n{n-k} +(m-l) \log \frac m{m-l}\r)(1+o(1))\\
&& +(n-k) \log \frac{ne}{n-k} + (m-l) \log \frac{me}{m-l}.
\end{eqnarray*}
Let us see that $(N-n)/n = N/n (1+o(1))$ and that
$$
(n-k)\log \left(\frac{n^2e}{(n-k)^2} \right) = n \left(1- \frac
kn\right)\l(1-2\log\l(1- \frac kn\r)\r) \leq \frac 2{\sqrt{e}}\, n,
$$
as $x(1-2\log(x))\leq 2/\sqrt{e}$ for all $x$ in $[0,1]$.


Let us denote $\mathcal{X}: = n\log(p^{-1})$ and $\mathcal{Y}: = m\log(q^{-1})$. We have
$$
\log\left(C_{N-n}^{n-k} \cdot C_{M-m}^{m-l}\cdot C_n^k \cdot
C_m^l\right) \leq \left(\left(1-\frac kn\right) \mathcal{X} +
\left(1-\frac lm\right) \mathcal{Y}+ \frac
2{\sqrt{e}}(n+m)\right)(1+o(1)).
$$
Analogously to \nref{alpha1} we have
\begin{eqnarray*}
\frac{(1-\delta_1)^2}2 a^2 nm \geq (1+\alpha_1) ( \cal{X}+\cal{Y} ),
\end{eqnarray*}
asymptotically.

Finally, we get, for large enough $n,\,m,\,N,\,M$
\begin{eqnarray*}
&&-\frac {(1-\delta_1)^2 a^2}2 (nm-kl) +\log\l(C_{N-n}^{n-k} \cdot C_{M-m}^{m-l}\cdot C_n^k \cdot C_m^l\r)\\
&\leq & - \alpha_1 \left( 1- \frac{kl}{nm}\right)(\mathcal{X}+\mathcal{Y})\\
&& -\left( 1- \frac{kl}{nm}\right)(\mathcal{X}+\mathcal{Y})+\l(\l(1-\frac kn\r)\mathcal{X} + \l(1-\frac lm\r) \mathcal{Y}+ \frac 2{\sqrt{e}}(n+m)\r)(1+o(1))\\
&\leq & -\frac{\alpha_1}2 \left( 1- \frac{kl}{nm}\right)(\mathcal{X}+\mathcal{Y})+\frac kn \l(\frac lm - 1\r)\mathcal{X} +
 \frac lm\l(\frac kn -1\r)\mathcal{Y}+ \frac 2{\sqrt{e}}(n+m)(1+o(1))\\
&\leq & -\frac{\alpha_1}2 \delta (\mathcal{X}+\mathcal{Y}) + \frac 2{\sqrt{e}}(n+m)(1+o(1)).
\end{eqnarray*}
Therefore, we replace this bound in \nref{eqT1} and get
\begin{eqnarray*}
&&\sum_{k=0}^n \sum_{l=0}^m I_{kl<(1-\delta)nm} T_{1,kl} \\
&\leq &  2\exp\left(
-\frac{\alpha_1}2 \delta (n\log(p^{-1})+m \log(q^{-1})) + \frac 2{\sqrt{e}}(n+m)(1+o(1))+\log (nm)
\right)=o(1).
\end{eqnarray*}

\bigskip

For $T_{2,kl}$, only the common elements of $C$ and $C_0$ play a role on the random variable $\sum_{C\cap C_0} \xi_{ij}$ and there are $C_n^k \cdot C_m^l$
 such choices. Note that we cannot have here neither $k=0$ nor $l=0$, as $T_{2,kl}=0$ in this cases. Therefore,
 \begin{eqnarray*}
\sum_{k=1}^n \sum_{l=1}^m I_{kl<(1-\delta)nm} T_{2,kl} &\leq &
\sum_{k=1}^n \sum_{l=1}^m C_n^k \cdot C_m^l P\left(Z >\frac{\delta_1
a(nm-kl)}{2 \sqrt{kl}}\right )
\\
& \leq &\sum_{k=1}^n \sum_{l=1}^m C_n^k \cdot C_m^l P\left(Z >\frac{\delta_1 \delta anm}{2 \sqrt{(1-\delta ) nm}} \right)\\
&\leq & \sum_{k=1}^n \sum_{l=1}^m \exp\left( -\frac{\delta_1^2
\delta^2 a^2 nm}{8(1-\delta)} + k \log\left(\frac {ne}k\right) + l
\log\left(\frac {me}l\right)
\right)\\
&\leq & \exp\left(-\frac{\delta_1^2 \delta^2 a^2 nm}{8(1-\delta)} +
 n  +  m  + \log(nm) \right) = o(1).
\end{eqnarray*}

Here, we have used the fact that $x\log(x^{-1})$ is bounded from
above by $e^{-1}$ for all $x \in [0,1]$ and used it for $x= k/n$ and
for $x=l/m$, respectively.
 Use \nref{alpha} in order to conclude.

\bigskip

Finally, for $T_{3,kl}$, we write that $-\sum_{C_0} \xi_{ij}/\sqrt{nm}$ behaves
like some standard Gaussian random variable $Z$ and get
\begin{eqnarray*}
\sum_{k=0}^n \sum_{l=0}^m T_{3,kl}
&\leq & \sum_{k=0}^n \sum_{l=0}^m \exp\left(- \frac{\delta_1^2 a^2 (nm-kl)^2}{8 nm} \right)
\\
& \le & \exp\left(- \frac{\delta_1^2 \delta^2 a^2 }{8 }nm +  \log
(nm) \right) = o(1),
\end{eqnarray*}
as $a ^2 {nm}$ tends to infinity faster than $\log(nm)$ due to \nref{alpha} in our setup.

In conclusion, the probability in \nref{weakdoublesum} tends to 0:
\begin{equation}\label{weak}
 P_{C_0} \left(\max_{k,\, l} \max_{kl<(1-\delta)nm}
\max_{C \in\, \CC_{nm,kl}}(\sum_{C\setminus  C_0} \xi_{ij} - \sum_{C_0 \setminus C} \xi_{ij} - a(nm-kl) ) > 0 \right) = o(1).
\end{equation}

\subsection{Large intersection}

Let us fix $k$ and $l$ such that $kl \geq (1-\delta)nm$, or, equivalently, $nm-kl \leq \delta nm$. Note that it implies both
$k \geq (1-\delta_1)n$ and $l \geq (1-\delta_1)m$ for some $\delta_1$ depending on $\delta$ small as $\delta \to 0$.
The case $n=k$ and $m=l$ gives an event with 0 probability.

We decompose as follows
\begin{eqnarray*}
&& \sum_{C\backslash C_0}\xi_{ij} - \sum_{C_0 \backslash C}\xi_{ij}
 =  \left(\sum_{(A\backslash A_0) \times B_0 }\xi_{ij} - \sum_{(A_0 \backslash A) \times B_0 }\xi_{ij} \right)\\
&&+ \left(\sum_{ A_0 \times (B \backslash B_0)}\xi_{ij}  - \sum_{ A_0 \times (B_0 \backslash B)}\xi_{ij}
\right) \\
&& + \left( \sum_{(A \backslash A_0) \times (B \backslash B_0)}\xi_{ij}   - \sum_{(A \backslash A_0)\times (B_0 \backslash B)}\xi_{ij}
+ \sum_{(A_0 \backslash A) \times (B_0 \backslash B)}\xi_{ij}   - \sum_{(A_0 \backslash A)\times (B \backslash B_0)}\xi_{ij}
\right)
\\
& = & S_1+S_2+S_3, \mbox{ say}.
\end{eqnarray*}
We shall bound from above as follows
\begin{eqnarray*}
&& P_{C_0} \left( \max_{k\geq (1-\delta_1)n} \max_{l\geq (1-\delta_1)m}
\max_{C \in\, \CC_{nm,kl}}(\sum_{C\setminus  C_0} \xi_{ij} - \sum_{C_0 \setminus C} \xi_{ij} - a(nm-kl) ) > 0\right)\\
&\leq & P_{C_0} \left(\max_{k\geq (1-\delta_1)n} \max_{l\geq (1-\delta_1)m} \max_{A \in\mathcal{C}_{n,k}} (S_1 - (1- {\tilde \delta})a(n-k) \frac{m+l}2 )>0\right)\\
&& +P_{C_0} \left(\max_{k\geq (1-\delta_1)n} \max_{l\geq (1-\delta_1)m}  \max_{B \in \mathcal{C}_{m,l}}(S_2 - (1-{\tilde \delta})a(m-l) \frac{n+k}2 ) >0 \right)\\
&& +  P_{C_0}\left(\max_{k\geq (1-\delta_1)n} \max_{l\geq (1-\delta_1)m} \max_{C \in\, \CC_{nm,kl}} (S_3 - {\tilde \delta} a(nm-kl) ) >0 \right),
\end{eqnarray*}
where $\mathcal{C}_{n,k}$ is the set of $n$ rows in $1,...,N$ having $k$ values in common with $A_0$ and similarly for $\mathcal{C}_{m,l}$ set
of $m$ columns in $1,...,M$ having $l$ values in common with $B_0$.
Moreover, the previous sum can be bounded from above by
\begin{eqnarray*}
&& \sum_{k\geq (1-\delta_1)n}P_{C_0} \left( \max_{A \in\mathcal{C}_{n,k}} S_1 > (1- {\tilde \delta})a(n-k) m(1-\delta_1/2) \right) \\
&&+\sum_{l\geq (1-\delta_1)m} P_{C_0} \left(  \max_{B \in \mathcal{C}_{m,l}}S_2 > (1-{\tilde \delta})a(m-l) n (1-\delta_1/2)  \right)\\
&& + \sum_{k\geq (1-\delta_1)n} \sum_{l\geq (1-\delta_1)m}P_{C_0}\left( \max_{C \in\, \CC_{nm,kl}}S_3 > {\tilde \delta} a(nm-kl) \right)\\
&=&
\sum_{k\geq (1-\delta_1)n} U_{1,k} + \sum_{l\geq (1-\delta_1)m}  U_{2,l}
+ \sum_{k\geq (1-\delta_1)n} \sum_{l\geq (1-\delta_1)m} U_{3,kl} \mbox{ say},
\end{eqnarray*}

Let us now deal with $ U_{1,kl} $. Note, first, that the case $k=n$ gives probability 0. For $(1-\delta_1/2) n \leq k \leq n-1$, we put $p_{n,N} = \sqrt{\log(N-n)}/(\sqrt{\log(N-n)}+\sqrt{\log (n)})$ and $q_{n,N} = 1-p_{n,N}$,
\begin{eqnarray*}
U_{1,k} &\leq &
P_{C_0} \left(\max_{A\in \mathcal{C}_{n,k}}\sum_{(A\backslash A_0) \times B_0 }\xi_{ij}  > (1- {\delta})(1-\delta_1/2)a(n-k)m p_{n,N}\right)\\
&& + P_{C_0} \left(\max_{A\in \mathcal{C}_{n,k}} \sum_{(A_0\backslash A) \times B_0 }(-\xi_{ij}) > (1- {\delta}) (1-\delta_1/2) a(n-k) m q_{n,N}\right)
\end{eqnarray*}
and, for some independent standard gaussian r.v. $Z_1$ and $Z_2$, using $l\geq (1-\delta_1)m$
\begin{eqnarray*}
U_{1,k}&\leq &
C_{N-n}^{n-k} P (Z_1  > (1- \delta) (1-\delta_1/2)p_{n,N} a \sqrt{(n-k)m})\\
&& + C_n^k P (Z_2  > (1- \delta) (1-\delta_1/2)q_{n,N} a \sqrt{(n-k)m})\\
&\leq & \exp\left( -\frac {(1-\tilde\delta)^2}2 \frac{a^2 m(n-k) \log(N-n)}{ (\sqrt{\log(N-n)}+\sqrt{\log (n)})^2}
+ \log(C_{N-n}^{n-k}) \right)\\
&& + \exp\left( -\frac {(1-\tilde\delta)^2}2 \frac{a^2 m(n-k) \log(n)}{ (\sqrt{\log(N-n)}+\sqrt{\log (n)})^2}
+ \log(C_{n}^{k}) \right),
\end{eqnarray*}
with $1-\tilde \delta = (1-\delta)(1-\delta_1/2)$. Note that $\log(C_{N-n}^{n-k}) \leq (n-k)\log(N-n)(1+o(1))$ and that $\log(C_{n}^{k}) \leq (n-k)\log(n)(1+o(1))$.
We obtain
\begin{eqnarray*}
U_{1,k} &\leq & \exp \left(-(n-k)\log(N-n)\left(
\frac {(1-\tilde\delta)^2}2 \frac{a^2 m}{ (\sqrt{\log(N-n)}+\sqrt{\log (n)})^2}
-(1+o(1)) \right) \right)\\
&& +
\exp \left(-(n-k)\log(n)\left(
\frac {(1-\tilde\delta)^2}2 \frac{a^2 m}{ (\sqrt{\log(N-n)}+\sqrt{\log (n)})^2}
-(1+o(1)) \right) \right).
\end{eqnarray*}
We use \nref{alpha}, for small enough $\delta$
\begin{eqnarray*}
(1-\tilde\delta)^2 a^2 m &\geq&
(1+2\alpha_2)2 (\sqrt{ \log(n)}+\sqrt{\log(N-n)})^2,
\end{eqnarray*}
for some $\alpha_2>0$ and this means
$$
\frac {(1-\tilde\delta)^2}2 \frac{a^2 m}{ (\sqrt{\log(N-n)}+\sqrt{\log (n)})^2}
-(1+o(1)) \geq 2\alpha_2 - o(1) \geq \alpha _2.
$$
Finally,
\begin{eqnarray*}
\sum_{ (1-\delta_1)n \leq k < n} U_{1,k } &\leq & \sum_{(1-\delta_1)n \leq k < n}(e^{- \alpha_2 \log(N-n) (n-k)}+e^{- \alpha_2 \log(n) (n-k)} )\\
& \leq & \sum_{1 \leq j \leq \delta_1 n}(e^{- \alpha_2 \log(N-n) j}+e^{- \alpha_2 \log(n) j} )\\
& = & (e^{- \alpha_2 \log(N-n)}+e^{- \alpha_2 \log(n)})(1+o(1)) =o(1).
\end{eqnarray*}

\bigskip

The term $U_{2,l}$ is similar.

\bigskip

As for the last term, $U_{3,kl}$, we compare each sum in $S_3$ to $\tilde{\delta} a(nm-kl)/4$.
The most difficult (the largest) upper bound is for the first sum, as it gives the largest number of choices $C_{N-n}^{n-k} C_{M-m}^{m-l}$.
Note that this term is 0 if $k=n$ or $l=m$. Therefore, we only explain this term, for $k\leq n-1$ and $l\leq m-1$,
\begin{eqnarray*}
U_{31,kl}&=&P_{C_0} \left(\max_{C \in\, \CC_{nm,kl}}\sum_{(A\backslash A_0)\times (B \backslash B_0)}\xi_{ij}
> \frac{\tilde\delta}{4} a(nm-kl) \right)\\
&\leq & C_{N-n}^{n-k} C_{M-m}^{m-l} \exp\left( -\frac {(\tilde \delta/4)^2} 2 \frac{a^2(nm-kl)^2}{(n-k)(m-l)}\right)\\
&\leq & \exp \left( -\frac{(\tilde \delta/4)^2 a^2 (n(m-l)P_{k,n}+(n-k)m P_{l,m})^2}{2(n-k)(m-l)} \right.\\
&& \left. +(n-k) \log(N-n) +(m-l) \log(M-m)\right),
\end{eqnarray*}
where $P_{k,n} = 1-(n-k)/(2n)$ and $P_{l,m} = 1-(m-l)/(2m)$. Recall that $n-k \leq \delta_1 n$ and that $m-l \leq \delta_1 m$. We get
\begin{eqnarray*}
U_{31,kl} & \leq &
\exp \left( -\frac{(\tilde \delta/4)^2 a^2}2 \left(H + 2nm P_{k,n} P_{l,m} \right) + \delta_1 (n\log(N-n) +m \log(M-m))\right),
\end{eqnarray*}
where $$
H=\frac{n^2}{n-k} (m-l) P^2_{k,n} +(n-k) \frac{m^2}{m-l} P^2_{l,n} \geq \frac 1{\delta_1} (n P^2_{k,n} + m P^2_{l,n}).
$$
Recall that $P_{k,n} \geq 1-\delta_1/2$ and $P_{l,m} \geq 1-\delta_1/2$. We get
for $(\tilde \delta/4)^2=\delta_1$:
\begin{eqnarray*}
U_{31,kl} & \leq &
\exp \left( -\frac{ a^2}2 (n P^2_{k,n} + m P^2_{l,n})- \delta_1(a^2nm P_{k,n} P_{l,m} - (n\log(N-n) +m \log(M-m)) \right),
\end{eqnarray*}
with
\begin{eqnarray*}
a^2nm P_{k,n} P_{l,m} &\geq &(1-\delta_1/2)^2 (\frac 12 a^2nm + \frac 12 a^2nm) \\
&\geq &(1-\delta_1/2)^2(1+\alpha) (n \log(n(N-n))+m \log(m(M-m))),
\end{eqnarray*}
by \nref{alpha}. By taking $\delta_1$ small enough, we may find $\delta_2>0$ such that $(1-\delta_1/2)^2(1+\alpha) \geq 1+\delta_2$.
This is enough to conclude that
$$
a^2nm P_{k,n} P_{l,m} - (n\log(N-n) +m \log(M-m)) >0
$$
and that
\begin{eqnarray*}
U_{31,kl} & \leq &
\exp \left( -\frac{ a^2}2 (n  + m )(1-\delta_1/2)^2  \right)\\
&\leq & \exp\left( -(1-\delta_1/2)^2(1+\alpha) (\log(m(M-m))+\log(n(N-n)))\right)\\
&\leq & \exp\left( -(1+\delta_2) (\log(m(M-m))+\log(n(N-n)))\right).
\end{eqnarray*}
In conclusion,
\begin{eqnarray*}
\sum_{(1-\delta_1)n\leq k <n} \sum_{(1-\delta_1) m\leq l < m} U_{31,kl} \leq \exp\left( -(1+\delta_2) \log((M-m)(N-n))-\delta_2 \log(nm) \right)=o(1).
\end{eqnarray*}

Here, we have proven that
\begin{equation}\label{strong}
 P_{C_0} \left( \max_{k,\, l} \max_{kl\geq (1-\delta)nm}
\max_{C \in\, \CC_{nm,kl}}(\sum_{C\setminus  C_0} \xi_{ij} - \sum_{C_0 \setminus C} \xi_{ij} - a(nm-kl) ) > 0 \right) = o(1).
\end{equation}
From \nref{strong} and \nref{weak} we deduce that the probability $P_{C_0}(\hat C^\star \ne C_0)$ tends to 0 and this concludes the proof of the upper bounds.
\endproof

\begin{remark}
We have investigated the upper limits of the selector $\hat C^\star  $
under the assumption that $s_{ij}=a, \ (i,j)\in C_0 $. It follows
that, when $s_{ij}\ge a, \ (i,j)\in C_0 $, statements of upper bounds
stated in this section are valid.

 Indeed, the random part of the expansion $Y_C-Y_{C_0} $
is independent of $s_{ij}$. The absolute value of the deterministic part
(the difference of expectations) attains its minimum when
$s_{ij}=a$.
\end{remark}

\section{Lower bounds}\label{lowb1}
Let \nref{lim} and \nref{Bsup}. We shall call the case when $B = A $
the case of severe sparsity, while the case where either $B = A_1 $
or $B = A_2$ will be designated by moderately sparse cases. Let us
first consider a set $\Theta$ of matrices having size $N \times M$
and containing $S_{C}$, for all $C \in \mathcal{C}_{nm}$, such that
$[S_C]_{ij} = a \cdot I((i,j) \in C)$. This set is on the border of
$\mathcal{S}_{nm,a}$, as we replace $[S_C]_{ij} \geq a$ with
equality, for all $(i,j) \in C$. The set $\Theta$ has $L = C_N^n
\cdot C_M^m$ elements. Let $P_0$ denote the likelihood of $N\times
M$ standard gaussian observations and, as previously,
 $P_C$ the likelihood of our observations under parameter $S_C$.
The minimax risk is bounded from below by the minimax risk over $\Theta$:
$$
\inf_{\hat C} \sup_{ S_{C} \in\,\CS_{nm,a} }P_{C}(\hat C(Y)\ne C) \geq \inf_{\hat C} \sup_{ S_{C} \in\,\Theta }P_{C}(\hat C(Y)\ne C) .
$$

\subsection{Severe sparsity}\label{sec:estimlb}
{\bf Proof of Theorem~\ref{estimlb} for severely sparse case}

In this case, we shall apply Theorem 2.4 in \cite{Tsyb}: if there exists $\tau >0$ and $0<\alpha < 1$ such that
$$
\frac 1L \sum_{S_C \in \Theta} P_C\left( \frac{dP_0}{dP_C} \geq \tau \right)\geq 1-\alpha,
$$
then
$$
\inf_{\hat C} \sup_{ S_{C} \in\,\Theta }P_{C}(\hat C(Y)\ne C) \geq \frac{\tau L}{1+ \tau L}(1-\alpha).
$$
In our model, the likelihood ratio is
\begin{equation}\label{liratio}
\frac{dP_0}{dP_C} = \exp \left( -a\sum_{(i,j) \in C} Y_{ij} +\frac{a^2 nm}2\right).
\end{equation}
This implies that
\begin{eqnarray*}
P_C\left( \frac{dP_0}{dP_C} \geq \tau \right) &=&
P_C \left(-a\sum_{(i,j) \in C} Y_{ij} +\frac{a^2 nm}2 \geq \log(\tau) \right)\\
&=& P_0\left(-\frac 1{\sqrt{nm}}\sum_{(i,j) \in C} \xi_{ij} -\frac{a\sqrt{nm}}2\geq \frac{\log(\tau)}{a \sqrt{nm}} \right)\\
&=& P\left(Z \geq \frac{\log(\tau)}{a \sqrt{nm}} + \frac{a\sqrt{nm}}2 \right),
\end{eqnarray*}
where $Z$ is standard gaussian. Let $z_{1-\alpha}$ be the quantile of probability $1-\alpha$ of a standard gaussian distribution, such that
 $P(Z\geq -z_{1-\alpha}) = 1-\alpha$. In order to check \nref{liratio}, we need $\log(\tau ) \leq -a^2 nm/2 - z_{1-\alpha} a\sqrt{nm}$.

On the one hand, if $a\sqrt{nm} =O(1)$ we take $\tau$ as solution of the equation $\log(\tau ) = -a^2 nm/2 - z_{1-\alpha} a\sqrt{nm}$.
 Therefore, we have $\tau \asymp 1$ and then
$$
\frac{\tau L}{1+ \tau L}(1-\alpha) \geq (1-\alpha)^2 >0, \quad \mbox{ as } L \to \infty.
$$

On the other hand, if $a\sqrt{nm} \to \infty$, we take $\tau^{-1} = L/\log(L)$, with $L=C_N^n C_M^m$, which gives $\tau L \to \infty$
and $\log(\tau^{-1}) \sim \log(L)$. We can prove that
$$
\log(\tau^{-1}) \geq \frac{a^2nm}2 +z_{1-\alpha} a\sqrt{nm} = \frac{a^2nm}2 \left(1+ \frac{2z_{1-\alpha}}{a\sqrt{nm}}\right).
$$
Indeed, we known that
$\log(L) \sim n\log(p^{-1})+m\log(q^{-1}) $ and, by assumption \nref{Bsup},
$$
\frac{a^2nm}{2(n\log(p^{-1})+m\log(q^{-1}))} \leq 1-\delta,
$$
asymptotically, for some $\delta >0$. It implies that
$$
\frac{a^2 nm}{2 \log(\tau^{-1})} \leq \left( 1+ \frac{2z_{1-\alpha}}{a\sqrt{nm}}\right)^{-1},
$$
asymptotically. This gives the lower bound
$$
\frac{\tau L}{1+ \tau L}  (1-\alpha) \geq (1-\alpha)^2 > 0.
$$
As $\alpha>0$ can be chosen arbitrarly small, we obtain the result
$$
\inf_{\hat C} \sup_{ S_{C} \in\,\Theta }P_{C}(\hat C(Y)\ne C) \to 1.
$$

\endproof

\subsection{Moderate sparsity}\label{sec:notcons}

\begin{lemma}\label{maxgauss}If $\eta_1,...,\eta_J$ are i.i.d. random variables with standard gaussian law, then
$$
\text{ if } t<1, \quad P(\max_{j=1,...,J} \eta_j \geq t\sqrt{2\log(J)}) \to 1, \mbox{ as } J\to \infty,
$$
and
$$
\text{ if } t>1, \quad P(\max_{j=1,...,J} \eta_j \geq t\sqrt{2\log(J)}) \to 0, \mbox{ as } J\to \infty.
$$
\end{lemma}

{\bf Proof} This Lemma is an obvious consequence of the limit behaviour of the normalized maximum of i.i.d. Gaussian random variables as follows:
$$
V_J:=\max_{j=1,...,J} \eta_j \sqrt{2 \log(J)} -2\log(J) +\frac 12 \log(\log(J))+\frac 14 \log(4\pi) \to^d U,
$$
where $U$ has the Gumbel law with distribution function $P(U\leq x )=\exp(-\exp(-x))$ for all real number $x$, see \cite{Emb}.
Therefore, if $t<1$,
\begin{eqnarray*}
P(\max_{j=1,...,J} \eta_j \geq t\sqrt{2\log(J)})  =
P(V_J \geq (t-1){2\log(J)}+\frac 12 \log(\log(J))+\frac 14 \log(4\pi) ),
\end{eqnarray*}
which tends to 1 when $J \to \infty$.
The other limit is obtained by a similar argument.
\endproof

\bigskip

{\bf Proof of Proposition~\ref{notcons}} Let us assume that $\limsup A_1 < 1$ and treat the other case similarly.
This means that $A_1 \leq 1-\alpha$, for some fixed $0<\alpha <1$. Equivalently, $a\sqrt{m} \leq (1-\alpha)(\sqrt{2 \log(n)}+\sqrt{2\log(N-n)})$.

In this case we shall reduce the set of matrices $C$ to those matrices having the same columns as $C_0$ and $n-1$ rows in common with $C_0$.
 Then we sum up each line over these columns and reduce the problem to the vector case. Thus,
\begin{eqnarray*}
P_{C_0}(\hat C^\star \not = C_0) &=& P_{C_0}(\max_{C \in \mathcal{C}_{nm}} \sum_C Y_{ij} -\sum_{ C_0}Y_{ij}>0) \\
&\geq & P_{C_0}(\max_{C =A\times B_0} \sum_C Y_{ij} -\sum_{ C_0}Y_{ij}>0) \\
&\geq & P_{C_0}(\max_{A} \sum_A Y_{i\cdot} -\sum_{ A_0}Y_{i\cdot}>0),
\end{eqnarray*}
where the maximum over $A$ is taken over all sets of $n$ rows having $n-1$ rows in common with $A_0$ and
$$
Y_{i\cdot}:= \sum_{j\in B_0} Y_{ij} = am I(i\in A_0) +\sum_{j\in B_0}\xi_{ij}.
$$
Denote by $\eta_i = m^{-1/2}\sum_{j\in B_0}\xi_{ij}$ for $i=1,...,N$, which are i.i.d. random variables of standard gaussian law. Therefore, we get
\begin{eqnarray*}
P_{C_0}(\hat C^\star \not = C_0)
& \geq & P_{C_0} (\max_{A} \sum_A \eta_i - \sum_{A_0}(\eta_i+a\sqrt{m})>0)\\
& \geq & P_{C_0} (\max_{i \not\in A_0} \eta_i + \max_{k\in A_0} (-\eta_k) >a\sqrt{m})\\
& \geq & P_{C_0} (\max_{i \not \in  A_0} \eta_i + \max_{k\in A_0} (-\eta_k) > (1-\alpha)(\sqrt{2\log(N-n)} + \sqrt{2 \log(n)}))\\
&=& 1- P_{C_0}(\max_{i \not \in  A_0} \eta_i + \max_{k\in A_0} (-\eta_k) \leq (1-\alpha)(\sqrt{2\log(N-n)} + \sqrt{2 \log(n)})),
\end{eqnarray*}
by the assumption on $A_1$.
Moreover
\begin{eqnarray*}
&&P_{C_0}(\max_{i \not \in  A_0} \eta_i + \max_{k\in A_0} (-\eta_k) \leq (1-\alpha)(\sqrt{2\log(N-n)} + \sqrt{2 \log(n)}))\\
&\leq & P_{C_0}(\max_{i \not \in  A_0} \eta_i \leq (1-\alpha)\sqrt{2\log(N-n)})
+ P_{C_0}(\max_{k\in A_0} (-\eta_k) \leq (1-\alpha) \sqrt{2 \log(n)} ),
\end{eqnarray*}
which tends to 0, by Lemma~\ref{maxgauss}.
\endproof


{\bf Proof of Theorem~\ref{estimlb} for moderately sparse case}.

In this case we check that the minimax risk is bounded from below by the risk of the maximum likelihood estimator $\hat C^\star$ and
 that its risk tends to 1 under our assumptions by Proposition~\ref{notcons}.
Let us see that
\begin{eqnarray*}
\inf_{\hat C} \sup_{ S_{C} \in\,\Theta }P_{C}(\hat C(Y)\ne C)
&\geq &\inf_{\hat C} \frac 1L
\sum_{k=1}^L P_{C_k}(\hat C(Y)\ne C_k)\\
&\geq & \inf_{\hat C} \left( 1-\frac 1L
\sum_{k=1}^L P_{C_k}(\hat C(Y)= C_k)\right)\\
&\geq & 1-\sup_{\hat C} \frac 1L
\sum_{k=1}^L E_{0}(I(\hat C(Y)= C_k)\frac{dP_{C_k}}{dP_0}(Y)),
\end{eqnarray*}
where $L=C_N^n C_M^m$ is the number of elements in $\Theta$.
In the previous supremum, we may replace the arbitrary measurable function $\hat C(Y)$ by a test function
$\psi(Y)$ taking values in $1,...,L$. The test maximising
$$
\sup_{\psi(Y)} \frac 1L
\sum_{k=1}^L E_{0}(I(\psi (Y)= k)\frac{dP_{C_k}}{dP_0}(Y))
$$
will choose $k$ such that $C_k$ has maximal likelihood:
$\{Y: \frac{dP_{C_k}}{dP_0}(Y)\geq \frac{dP_{C_j}}{dP_0}(Y), \mbox{ for all } j =1,...,L\}$.
 Thus, we get the risk of a maximum likelihood estimator,
\begin{eqnarray*}
\inf_{\hat C} \sup_{ S_{C} \in\,\Theta }P_{C}(\hat C(Y)\ne C)
&\geq & 1-\frac 1L \sum_{k=1}^L P_{C_k}(\hat C^\star (Y) = C_k)\\
&\geq & \frac 1L \sum_{k=1}^L P_{C_k}(\hat C^\star (Y) \ne C_k),
\end{eqnarray*}
which tends to 1 by Proposition~\ref{notcons}.

\endproof



\begin{thebibliography}{99}



\bibitem{ACCD11} Arias-Castro, E.,  Cand\`es, E.J. and Durand, A.
(2011) Detection of an anomalous clusters in a network. \textit{Ann.
Statist.} {\bf 39} (1), 278-304.


\bibitem{arias} 
Arias-Castro, E., Cand\`es, E. J. and Plan, Y. (2010) Global Testing
and Sparse Alternatives: ANOVA, Multiple Comparisons and the Higher
Criticism. {\texttt arXiv:1007.1434}.

\bibitem{ACDH}
Arias-Castro, E., Donoho, D. L. and Huo, X. (2005) Near-optimal
detection of geometric objects by fast multiscale methods. {\it IEEE
Transactions on Information Theory} {\bf 51}, 7, 2402-2425.

\bibitem{ACL}
Arias-Castro, E., Lounici, K. (2012)
Variable selection with exponential weights and $\ell_0$-penalization.
arxiv:1208.2635.

\bibitem{BH}
Benjamini, Y. and Hochberg, Y. (1995) Controlling the false discovery rate: a practical and powerful approach to multiple testing. {\it J. Roy. Statist. Soc. Ser. B}, {\bf 57}, 289-3300.

\bibitem{BRT}
Bickel, P.J., Ritov, Y. and Tsybakov, A.B.  (2009) Simultaneous
analysis of Lasso and Dantzig selector. {\it Annals of Statistics},
{\bf 37}, 4, 1705-1732.

\bibitem{ButIn}
Butucea, C., Ingster, Yu. I. (2012). Detection of a sparse submatrix
of a high-dimensionale noisy matrix. \textit{Bernoulli}.

\bibitem{BG}
Butucea, C. and Gayraud, G. (2013). Sharp detection of smooth signals in a high-dimensional sparse matrix with indirect observations. arxiv:1301.4660

\bibitem{CJL} Cai, T., Jin, J. and Low, M. (2007) Estimation and confidence sets for sparse normal mixtures. \textit{Ann. Statist.}, {\bf 35}, 2421-2449.

\bibitem{DJ04} Donoho, D. and Jin, J. (2004) Higher criticism for detecting
sparse heterogeneous mixtures. \textit{Ann. Statist.} {\bf 32},
962-994.

\bibitem{Emb} Embrechts, P., Kl\"uppelberg, C and Mikosch, T. (1997)
Modelling extremal events: for insurance and finance.
{\it Springer}.



\bibitem{I97} Ingster, Yu.I. (1997)
Some problems of hypothesis testing leading to infinitely divisible
distributions. \textit{Math. Methods of Stat.} {\bf 6}, 47-69.

\bibitem{ingst} Ingster, Yu.I. and Stepanova, N.A. (2012)
Adaptive selection of sparse regression function components.
\textit{Zapiski Nauchn. Sem. POMI} {\bf ZAI} (in Russian).

\bibitem{IS02b} Ingster, Yu.I. and Suslina, I.A. (2002)
On a detection of a signal of known shape in multichannel system.
\textit{Zapiski Nauchn. Sem. POMI} {\bf 294},  88-112 (in Russian,
Transl. \textit{J. Math. Sci.} {\bf 127}, 1723-1736).



\bibitem{KLT}
Koltchinskii,V., Lounici,K. and Tsybakov, A.B. (2011) Nuclear norm
penalization and optimal rates for noisy low rank matrix completion.
{\it Annals of Statistics, to appear.}

\bibitem{SN}
Sun, X. and Nobel, A.B. (2010) On the maximal size of Large-Average
and ANOVA-fit Submatrices in a Gaussian Random Matrix.
\textit{ArXiv: 1009.0562v1}

\bibitem{SWPN}
Shabalin, A.A., Weigman, V.J., Perou, C.M. and Nobel, A.B. (2009).
Finding Large Average Submatrices in High Dimensional Data.
\textit{Annals of Applied Statistics} 3, 985-1012.

\bibitem{Tsyb}
Tsybakov, A.B. (2009) Introduction to nonparametric statistics. {\it Springer Series in Statistics},
Springer, New-York.

\bibitem{Verzelen}
Verzelen, N. (2012). Minimax risks for sparse regressions:
Ultra-high dimensional phenomenons.
{\it Electron. J. Stat.} 6, 3890.

\end{thebibliography}
\end{document}